\documentclass[12pt,centertags,oneside]{amsart}

\usepackage[utf8]{inputenc}  
\usepackage[T1]{fontenc} 

\usepackage[backend=bibtex, style=alphabetic]{biblatex}
\renewbibmacro{in:}{}

\usepackage{typearea}
\usepackage[hidelinks]{hyperref}
\usepackage{stmaryrd}
\usepackage{amssymb}
\usepackage{a4wide}
\usepackage[mathscr]{eucal}
\usepackage{mathrsfs}
\usepackage{dsfont} 
\usepackage{enumitem} 
\usepackage{typearea}
\usepackage{charter}
\usepackage{pdfsync}
\usepackage[a4paper,width=16.2cm,top=3cm,bottom=3cm]{geometry}
\usepackage{amsmath}
\usepackage{amstext}
\usepackage{amsthm}
\usepackage{amscd}

\usepackage{color}
\usepackage{xcolor}

\numberwithin{equation}{section}

\newtheorem{theorem}{Theorem}[section]
\newtheorem{theo}[theorem]{Theorem}
\newtheorem{proposition}[theorem]{Proposition}

\theoremstyle{definition}
\newtheorem{definition}[theorem]{Definition}

\newcommand{\cali}[1]{\mathscr{#1}}

\newcommand{\supp}{{\rm supp}}

\newcommand{\id}{{\rm id}}

\newcommand{\Jac}{{\rm Jac}}

\newcommand{\C}{\mathbb{C}}

\newcommand{\N}{\mathbb{N}}

\renewcommand{\P}{\mathbb{P}^k}

\renewcommand{\epsilon}{\varepsilon}

\title{On the stability of holomorphic families of endomorphisms of $\P$}

\author{Fran\c cois Berteloot and Xavier Buff}
\address{Institut de Mathématiques de Toulouse; UMR5219; Université de Toulouse; CNRS; UPS, F-31062 Toulouse Cedex 9, France}
\email{francois.berteloot@math.univ-toulouse.fr}
\email{xavier.buff@math.univ-toulouse.fr}

\date{}

\begin{document}

\begin{abstract}
In the context of holomorphic families of $\P$ endomorphisms, we show that various notions of stability are equivalent. This allows us to both extend and simplify the architecture of the proof of certain results of \cite{BBD}.

\end{abstract}

\maketitle



\section{Introduction and results}
A \emph{holomorphic family of degree $d\ge 2$ endomorphisms on $\P$}, parametrized by a complex manifold $M$, is a holomorphic map $f:M\times \P\to M\times \P$, of the form $(\lambda,z)\mapsto (\lambda,f_\lambda(z))$ such that $d\ge 2$ is the common algebraic degree of the maps $f_\lambda$ as endomorphisms of $\P$.

For such families, the simplest notion of dynamical stability concerns repelling cycles and is referred to as  \emph{weak stability}. We are going to define various versions of it.
For every integer $n$ and every parameter $\lambda$, we denote by  ${\cali R}_n(\lambda)$
the set of $n$-periodic repelling points of $f_\lambda$ which belong to the Julia set of $f_\lambda$. 
Let us recall that, by the Briend-Duval equidistribution theorem \cite{BrDu}, $\vert{ \cali R}_n(\lambda)\vert \sim d^{kn}$ for every $\lambda \in M$.

\begin{definition}
A family $f$ is said  \emph{partially weakly stable}  if there exists a sequence of trivial holomorphic laminations
$({\cali L}^r_n)_{n\in\N}$ in $M\times \P$ such that  ${\cali L}^r_n \cap (\{\lambda\}\times \P)\subset  {\cali R}_n(\lambda)$ for every $\lambda\in M$, $f({\cali L}^r_n )\subset {\cali L}^r_n $
and $\limsup _n {\vert{\cali L}^r_n\vert}{d^{-kn}}>0$.
Such a family is said \emph{asymptotically weakly stable} if moreover $\lim _n {\vert{\cali L}^r_n\vert}{d^{-kn}} =1$, and  \emph{weakly stable} if
${\cali L}^r_n \cap (\{\lambda\}\times \P)= {\cali R}_n(\lambda)$ 
 for every $\lambda \in M$.
\end{definition}

Roughly speaking, weak stability means that all repelling cycles move holomorphically.
In dimension $k=1$, Ma\~{n}\'e-Sad-Sullivan \cite{MSS} and, independently, Lyubich \cite{Lyu} have shown
that weak stability is equivalent to the existence of a holomorphic motion of the full Julia sets. One of their main tool is the so-called $\lambda$-lemma. Combining this lemma with the ergodicity of the equilibrium measure of $f_\lambda$,
one may see that partial weak stability is actually sufficient \cite{Ber}.

In higher dimensions, the classical approach fails and new tools have to be introduced.
This has been done in \cite{BBD} where, in particular, a probabilistic notion of stability called $\mu$-stability has been defined (see Section {\ref{SecEW}).
 It is proved in \cite{BBD} that the notions of weak stability
and $\mu$-stability are equivalent within the full family ${\cali H}_d(\P)$ of degree $d$ holomorphic  endomorphisms of $\P$, or within
arbitrary families of endomorphisms of $\mathbb{P}^{2}$. Using further techniques, Bianchi \cite{Bia} proved the equivalence of 
asymptotic weak stability and $\mu$-stability in any dimension and in the more general setting of holomorphic  families of polynomial like mappings with large topological degree.\\

Our aim in this note is to prove the equivalence of the notions of weak stability and $\mu$-stability. Our main results are the following.
In section 3 we prove the
\begin{proposition}\label{PropMT}
Let  $M$ be  a connected  complex manifold and $f:M\times \P \to M\times \P$ be  a degree $d\ge 2$ holomorphic family  of endomorphisms of $\P$.
If $f$ is $\mu$-stable then $f$ is weakly stable. 
\end{proposition}

Combining the above result with some of \cite{BBD} and \cite{BB}, we then show in Section 4 that all the notions of stability are equivalent  when the parameter space is simply connected (see Theorem \ref{theoStab} for a complete set of equivalences).

\begin{theorem}\label{ThmEquiv}
Let  $M$ be  a simply connected  complex manifold and $f:M\times \P \to M\times \P$ be  a degree $d\ge 2$ holomorphic family  of endomorphisms of $\P$.
 
Then: $f$ is $\mu$-stable  $\Leftrightarrow  f$ is weakly stable  $\Leftrightarrow  f$ is partially weakly stable.

\end{theorem}

Finally, we should mention that the above result also covers the fact that the bifurcation of repelling cycles is a  chain reaction in the sense of \cite{Ber}.

\section{Equilibrium webs, laminations, and $\mu$-stability}\label{SecEW}
In all this section $f$ is a given degree $d$ holomorphic family of endomorphisms of $\P$.
The critical set of $f$ is denoted $C_f$ and, for each parameter $ \lambda \in M$, the equilibrium measure of $f_\lambda$ is denoted $\mu_\lambda$
while its support, the (small) Julia set of $f_\lambda$, is denoted $J_\lambda$.\\

We endow $\P$ with the spherical distance $d_{\P}$ and the space ${\cali O}(M,\P)$ of holomorphic maps from $M$ to $\P$
 with the distance $d_{uloc}$ of local uniform convergence, which makes  $({\cali O}(M,\P),d_{uloc})$ a complete separable metric space. The graph of $\gamma \in {\cali O}(M,\P)$ is denoted $\Gamma_\gamma$.
The holomorphic map $f$ clearly induces a  continuous selfmap $${\cali F} :{\cali O}(M,\P) \to {\cali O}(M,\P)$$  which is defined by 
${\cali F}\cdot \gamma (\lambda):= f_\lambda(\gamma(\lambda))$ for any $\gamma \in {\cali O}(M,\P)$ and $ \lambda \in M.$
We will also use  the evaluation maps $$e_\lambda : {\cali O}(M,\P) \to \P$$ defined by 
$e_\lambda(\gamma):=\gamma(\lambda)$ for any $\gamma \in {\cali O}(M,\P)$ and $\lambda \in M$. \\

An important role is played by the two  following  (possibly empty)  closed $\cali F$-invariant subspaces of ${\cali O}(M,\P)$
$${\cali J}:=\{\gamma\in {\cali O}(M,\P)\;\colon\: \gamma(\lambda) \in J_\lambda,\; \forall \lambda \in M   \},$$
$$\cali J_s := \{\gamma \in \cali J \: : \: \Gamma_\gamma \cap (\cup_{m\ge 0} f^{-m}(\cup_{n\ge 0}f^n(C_f))) \ne \emptyset \}.$$

As we shall see, the $\mu$-stability of the family $f$
amounts to say  that there exists an ergodic dynamical system of the form $({\cali J},{\cali M}, {\cali F})$
which does not interact with the critical dynamics of $f$.

\begin{definition} An \emph {equilibrium web} for the family $f:M\times\P\to M\times \P$ is a probability measure 
$\cali M$ on
${\cali O}(M,\P)$ such that:
\begin{itemize}
\item[1.] ${\cali F}_* {\cali M} =\cali M$,
\item[2.] $\supp \cali M$ is compact in $({\cali O}(M,\P),d_{uloc})$,
\item[3.] $(e_{\lambda})_* \cali M=\mu_\lambda, \; \forall \lambda \in M$.
\end{itemize}

An equilibrium web $\cali M$ is said  \emph{acritical} if $\cali M(\cali J_s) = 0$, and \emph{ergodic} if the
dynamical system $({\cali O}(M,\P),{\cali F},{\cali M})$ is ergodic.

An \emph {equilibrium lamination} for the family $f:M\times\P\to M\times \P$ is a relatively compact subset $\cali L$  of $\cali J$ such that
\begin{itemize}
\item[1.] $\Gamma_\gamma \cap \Gamma_{\gamma'} =\emptyset$ for every distinct $\gamma,\gamma' \in {\cali L}$,
\item[2.] $\mu_\lambda\{\gamma(\lambda)\;\colon\;\gamma \in {\cali L}\}=1$ for every $\lambda\in M$,
\item[3.]    $\Gamma_\gamma$ does not meet the grand orbit of the critical set of $f$ for every    $\gamma\in {\cali L}$,
\item[4.] the map ${\cali F}:{\cali L}\to {\cali L}$ is $d^k$ to $1$. 
\end{itemize}

An equilibrium lamination is said \emph{subordinated  to an equilibrium web} $\cali M$ if 
${\cali M}({\cali L})=1$.

 A family $f$ is said to be \emph{$\mu$-stable} if it admits an  equilibrium lamination. 
\end{definition}

 It turns out that every equilibrium web is supported in $\cali J$ 
 (see the comment after Proposition 2.3 in \cite{BB}). Let us also mention that Bianchi and Rakhimov \cite {BR} have recently shown that the
 $\mu$-stability, i.e. the stability of the maximal entropy measure,  implies a similar property for all measures with entropy strictly bigger than $(k-1)\ln d$.\\
 
 The interplay between equilibrium webs and laminations is given by the following fundamental result from \cite{BBD}.
 
 \begin{theorem}\label{TheoInter}
 Let  $M$ be  a simply connected  complex manifold and $f:M\times \P \to M\times \P$ be  a degree $d\ge 2$ holomorphic family  of endomorphisms of $\P$. 
 \begin{itemize}
 \item[1)] If $f$ admits an acritical and ergodic equilibrium web $\cali M$ then there exists an equilibrium lamination 
${\cali L}$ for $f$ which is subordinated to ${\cali M}$, moreover ${\cali M}({\cali L} \Delta {\cali L'})=0$ for any other
equilibrium lamination ${\cali L}'$ of $f$,
\item[2)] if  $f$ admits an  equilibrium lamination $\cali L$ then there exists an acritical and ergodic equilibrium web $\cali M$ of $f$  to which
${\cali L}$ is subordinated.
 \end{itemize}
 \end{theorem}
 
 The first assertion is Theorem 4.1 of \cite{BBD}, the second one is only implicit there and we thus prove it below.
  \proof
 Pick $\gamma_0 \in {\cali L}$ and, for every $n\in \N$, set ${\cali M}_n:=\frac{1}{n}\sum_{i=1}^n d^{-ki} \sum_{{\cali F}^{\circ i} \sigma =\gamma_0} \delta_\sigma$. Then any weak limit $\cali M$ of $({\cali M}_n)_n$ is an equilibrium web for $f$
 (see \cite{BBD} Proposition 2.2). We will show that $\cali M$ yields an acritical and ergodic equilibrium web for $f$ to which $\cali L$ is subordinated.
Let us first check that for every $k\in \N$ and every $\gamma \in \supp {\cali M}$ one has:
$$\Gamma_\gamma \cap f^{\circ k}(C_f)\ne \emptyset \Rightarrow \Gamma_\gamma \subset f^{\circ k}(C_f).$$
Indeed, if this were not the case, by Hurwitz theorem, we could find some $n\in \N$ and some   $\sigma \in \supp {\cali M_n}$   such that   $\Gamma_\sigma \cap  f^{\circ k}(C_f)\ne \emptyset$, and therefore
$f^{\circ  (i+k)} (C_f)\cap \Gamma_{\gamma_0} \ne \emptyset$ for some $i\le n$ which, as $\gamma_0\in {\cali L}$,  is impossible. 

Now, for any fixed $\lambda_0\in M$ we get
\begin{eqnarray*}
&&{\cali M}(\{\gamma\in {\cali J}\;\colon \; \Gamma_\gamma \cap (\cup_{k\ge 0} f^{\circ k}(C_f))\ne \emptyset\})=
{\cali M}(\{\gamma\in {\cali J}\;\colon \; \Gamma_\gamma \subset (\cup_{k\ge 0} f^{\circ k}(C_f))\})
\le\\
&& {\cali M}(\{\gamma\in {\cali J}\;\colon \; (\lambda_0,\gamma(\lambda_0))\in (\cup_{k\ge 0} f^{\circ k}(C_f))\})
= (e_{\lambda_0})_* {\cali M}\cup_{k\ge 0} f_{\lambda_0}^{\circ k}(C_{f_{\lambda_0}})=0 
\end{eqnarray*}
where the last equality follows from the fact that $ (e_{\lambda_0})_* {\cali M}=\mu_{\lambda_0}$ 
gives no mass to pluripolar subsets of $\P$. Then ${\cali M}({\cali J}_s)=0$ follows from the $\cali F$-invariance of $\cali M$ and therefore, according to the  Proposition 2.4 of \cite {BBD}, $\cali M$ can be replaced by an acritical  and ergodic  equilibrium web, which we still note $\cali M$.

By the first assertion of Theorem \ref{TheoInter}, there exists an equilibrium lamination ${\cali L}'$ of $f$ such that
${\cali M}({\cali L}')=1$ and ${\cali M}({\cali L}\Delta {\cali L}')=0$. It follows that ${\cali M}({\cali L})=1$.\qed\\

\section{Proof of Proposition \ref{PropMT}}

Let $\cali L$ be an equilibrium lamination, subordinated to an equilibrium web $\cali M$, for $f$.
Let $\lambda_0 \in M$ and $z_0\in J_{\lambda_0}$ be such that $z_0$ is $p$-periodic and repelling for $f_{\lambda_0}$. Without any loss of generality we may assume that $p=1$. According to \cite{BBD} Lemma 2.5, there exists $\gamma_0\in {\cali J}$ such that $\gamma_0(\lambda_0)=z_0$, ${\cali F}(\gamma_0)=\gamma_0$ and   $\gamma_0(\lambda)$ is a fixed repelling point of $f_\lambda$ for $\lambda$ sufficiently close to $\lambda_0$.
We must show that $\gamma_0(\lambda)$ is repelling for every $\lambda\in M$.

We proceed by contradiction. If this is not the case, then there exists $\lambda_1 \in M$ such that at least one of the eigenvalues of $f'_{\lambda_1} (\gamma_0(\lambda_1))$ has a  modulus smaller or equal to $1$. 

The fixed point $\gamma_0(\lambda
_1)$  is non repelling, 
let us first show that after possibly moving 
$\lambda_1$ a little, it  becomes hyperbolic.
Let $Z:=\{(\lambda, w)\in M\times \C\;\colon\; \det \left(f'_\lambda (\gamma_0(\lambda)) - w\;\id\right)=0\}$. The canonical projection $\pi_M:Z\to M$
has degree $d$  for some $d\le k$ and $\min_{\pi_M^{-1}(\lambda_1)} \vert w\vert \le 1$. If $\min_{\pi_M^{-1}(\lambda_1)} \vert w\vert = 1$, we cannot have $\min_{\pi_M^{-1}(\lambda)} \vert w\vert \ge  1$ on some neighbourhood of $\lambda_1$ since otherwise
 the maximum modulus principle, applied to the restriction of the function $\frac{1}{w}$ to $Z$, would imply that $\min_{\pi_M^{-1}(\lambda_0)} \vert w\vert \le  1$. We may thus slightly move $\lambda_1$ so that  $\min_{\pi_M^{-1}(\lambda_1)} \vert w\vert <  1$ and  $\pi_M$ has exactly $d$ distinct preimages at $\lambda_1$. Now,  $\pi_M^{-1} (\lambda)$ is given by the graphs of $d$  holomorphic functions $w_j(\lambda)$ above a small neighbourhood of $\lambda_1$ and, since $\min_{\pi_M^{-1}(\lambda_0)} \vert w\vert >1$, any function $w_j$ is not constant if $\vert w_j(\lambda_1)\vert =1$.
Thus, after possibly moving $\lambda_1$ again, $\gamma_0(\lambda_1)$ stays non repelling and becomes hyperbolic.

Let  $W$ be  the local unstable manifold of $f_{\lambda_1}$ at $\gamma_0(\lambda_1)$. It is  a proper analytic subset  of some neighbourhood of $\gamma_0(\lambda_1)$ in $\P$.
Let $\Omega_0:=\{z\in \P\;\colon\; d_{\P}(z,z_0) <\alpha_0\}$ be a neighbourhood of $z_0$ in $\P$ and ${\cali L}_0$, ${\cali L}_1$ be the following subsets of $\cali L$
$${\cali L}_0:=\{\gamma\in \cali L\;\colon\; \gamma(\lambda_0) \in \Omega_0\}=e_{\lambda_0}^{-1}(\Omega_0)\cap {\cali L}$$
$${\cali L}_1:=\{\gamma\in \cali L\;\colon\; \exists n\in \N \;\textrm{such that}\; \gamma(\lambda_1) \in f_{\lambda_1}^{\circ n}(W)\}=\displaystyle\bigcup_{n\in \N}
\; e_{\lambda_1}^{-1}(f_{\lambda_1}^{\circ n}(W)) \cap {\cali L}.$$
Let us show  that
$${\cali M}({\cali L}_0)>0\;\textrm{and}\;  {\cali M}({\cali L}_1)=0.$$
As ${\cali M}({\cali L})=1$ we have ${\cali M}({\cali L}_0)={\cali M}(e_{\lambda_0}^{-1}(\Omega_0))$
and  ${\cali M}({\cali L}_1) \le \sum_{n\in \N} {\cali M}(  e_{\lambda_1}^{-1}(f_{\lambda_1}^{\circ n}(W))$. Since $z_0\in J_{\lambda_0}$ we have ${\cali M}({\cali L}_0)=(e_{\lambda_0})_* {\cali M}(\Omega_0)=\mu_{\lambda_0}(\Omega_0)>0$. Since $\mu_{\lambda_1}$ gives no mass to pluripolar sets, ${\cali M}(e_{\lambda_1}^{-1}(f_{\lambda_1}^{\circ n}(W)) )=(e_{\lambda_1})_* {\cali M}(f_{\lambda_1}^{\circ n}(W))=\mu_{\lambda_1}(f_{\lambda_1}^{\circ n}(W))=0$ for every $n\in \N$ and thus ${\cali M}({\cali L}_1)=0$. 

Now, the expected contradiction will be obtained by 
showing that, for $\Omega_0$ suitably small, $${\cali L}_0\subset {\cali L}_1.$$

As $f_{\lambda_0}$ is repelling at $z_0$, there exists $0<a<1$, $r_0>0$, a neighbourhood $V_0$ of $\lambda_0$ in $M$, and an inverse branch $\varphi : T_0 \to T_0$  of $f$ which is defined on the neighbourhood 
$T_0:=\{(\lambda,z)\in V_0\times \P\;\colon\; d_{\P}(z,\gamma_0(\lambda)) <r_0\}$  of $(\lambda_0,z_0)$, such that
\begin{eqnarray}\label{InvBr}
d_{\P}(\gamma_0(\lambda), \varphi(\lambda,z)) \le a \; d_{\P}(\gamma_0(\lambda), z);\; \forall (\lambda,z)\in T_0.
\end{eqnarray}

Owing to the equicontinuity of the family $\cali L$, we may take the neighbourhood $V_0$ of $\lambda_0$  and the neighbourhood  $\Omega_0$ of $z_0$ small enough so that

\begin{eqnarray}\label{EqiL}
\forall \gamma \in {\cali L} ,\; \forall \lambda\in V_0\;\colon\; \gamma(\lambda_0) \in \Omega_0 \Rightarrow  (\lambda,\gamma(\lambda))\in T_0.
\end{eqnarray}

Let  $\gamma \in {\cali L}_0$. It follows from (\ref{EqiL}) and the definition of ${\cali L}_0$ that the map $\lambda \mapsto \varphi(\lambda,\gamma(\lambda))$ is well defined 
on $V_0$. 
Now, since the map ${\cali F} : {\cali L} \to {\cali L}$ is $d^k$ to $1$, this implies that  there exists 
$\gamma_{-1} \in \cali L$ such that $ \varphi(\lambda,\gamma(\lambda))=\gamma_{-1}(\lambda)$ on $V_0$. By (\ref{InvBr}), one sees that $\gamma_{-1}$ actually belongs to
${\cali L}_0$. We may thus iterate this construction and find a sequence $(\gamma_{-n})_n$ in ${\cali L}_0$ such that 
${\cali F}^{\circ n}(\gamma_{-n})=\gamma$ and $\lim_n \gamma_{-n}(\lambda)=\gamma_0(\lambda)$ for every $\lambda \in V_0$. By analyticity, and equicontinuity of $\cali L$, this implies that $\lim_n \gamma_{-n}(\lambda)=\gamma_0(\lambda)$ for every $\lambda \in M$ and, in particular, that $\lim_n \gamma_{-n}(\lambda_1)=\gamma_0(\lambda_1)$. This is only possible if $\gamma_{-n}(\lambda_1) \in W$ for $n$ big enough,
and thus $\gamma\in {\cali L}_1$.
We have shown that ${\cali L}_0\subset {\cali L}_1.$\qed\\

\section{Proof of Theorem \ref{ThmEquiv}}

We shall actually present  the following more complete result which combines the above Proposition \ref{PropMT}
with results of \cite{BBD} and {\cite{BB}.

\begin{theo}\label{theoStab}
Let $f:M\times \P\to M\times \P$  be a holomorphic family of endomorphisms of degree $d\ge 2$
on $\P$ parametrized by a simply connected complex manifold $M$. Then the following assertions are equivalent:
\begin{itemize}
\item[1)] $L$ is pluriharmonic on $M$,
\item[2)] $\Vert f^{\circ n}_* [C_f]\Vert_{U\times \P}= O(d^{(k-1)n})$
for every $U \Subset M$,
\item[3)] the ramification current $R_f :=\sum_{n\ge 0} d^{-kn} (f^{\circ n})_* [ f(C_f)]$
converges on $M\times \P$,
\item[4)] $f$ is {$\mu$}-stable,
\item[5)] $f$ is  weakly stable,
\item[6)] $f$ is  asymptotically weakly stable,
\item[7)] $f$ is  partially weakly stable.
\end{itemize}
\end{theo}

Let us recall that $[C_f]$ denotes  the current of integration on the critical set $C_f$ of $f$, and that the Lyapunov function  $L(\lambda):=\int_{\P} \ln \vert \Jac f_\lambda \vert\;\mu_\lambda$ is $p.s.h$ on $M$ (\cite{BaBe}) and coincides  with  the sum of the Lyapunov exponents of the ergodic dynamical system
$(J_\lambda,f_\lambda, \mu_\lambda)$.

One of the main threads in the proof of the Theorem \ref{theoStab}
is that the stability properties of $f$ are encoded by the function $L$, from which they can be extracted using the following two formulas
\begin{itemize}[label={}, leftmargin=*]
    \item  $\boldsymbol{dd^c L}$\textbf{-formula}:
    $dd^cL=\pi_{M*}\left((dd^c_{\lambda,z} g(\lambda,z) +\omega_{\P})^k\wedge[C_f]\right)$,
    \item \textbf{Approximation formula}:
    $L(\lambda)=\lim_n d^{-kn}\sum_{z\in{\cali R}_n(\lambda)} \ln \vert {\Jac} f_\lambda(z)\vert$.
\end{itemize}~

The first formula has been obtained by Bassanelli and Berteloot in \cite{BaBe},
it generalizes similar formulas in dimension one due to Przytycki \cite{Prz} and Manning \cite{Man}  for polynomials 
and DeMarco  \cite{DMa} for rational functions. It might be useful to stress that $g(\lambda,\cdot)$ is the  Green function of $f_\lambda$ and that $\mu_\lambda= (dd^c_{z} g(\lambda,z) +\omega_{\P})^k$ (see \cite{DS} page 176).

The second formula was proved by Berteloot, Dupont and Molino in \cite{BDM} and a simplified proof,
avoiding difficulties due to the possible resonances between the Lyapunov exponents, has been given by Berteloot and Dupont in \cite{BD}.\\

Let us now enter into details. The implications $2)\Rightarrow 3)$ and $5)\Rightarrow 6) \Rightarrow 7)$ are obvious.
Note that the convergence of the positive current $R_f$ means that every point in $M\times \P$ admits a neighbourhood $U$ such that the series $\sum_{n\ge 0} \Vert \mathds{1}_U d^{-kn} (f^{\circ n})_* [ f(C_f)]\Vert$ converges.

$1) \Rightarrow 2)$ follows immediately from the following estimate which  is a direct consequence of the $dd^c L$-formula
(\cite[Lemma 3.13]{BBD} Lemma 3.13).
There exists a positive constant $\alpha$, only depending on $k$ and $\dim_\C M$,
such that
$\Vert f^n_* [C_f]\Vert_{U\times \P}= \alpha d^{kn} \Vert dd^c L\Vert_U + O(d^{(k-1)n})$
for every relatively compact open subset $U$ of $M$.

$3) \Rightarrow  4)$ This is based on the key result of \cite{BB}. Namely, the convergence of the ramification current implies that any  $\lambda_0\in M$ has a neighbourhood $D_0$
such that the restricted family $f\vert_{D_0\times \P}$ admits an acritical equilibrium web $\cali M$ (\cite{BB} Theorem 1.4 and Lemma 2.4).
Then, by the Proposition 2.4 in \cite{BBD} this web can be assumed to be ergodic, and thus $f\vert_{D_0\times \P}$ admits an equilibrium lamination (Theorem \ref{TheoInter}.)
As $M$ is simply connected,
$f$ itself admits an equilibrium lamination.

 $4)\Rightarrow 5)$ is Proposition \ref{PropMT}.
 
$7)\Rightarrow 1)$ This has been proved in \cite{Ber}, we reproduce here  the argument.   By Banach-Alaoglu theorem, there exists a subsequence $(n_q)_q$ and a compactly supported positive  measure $\cali M$ of mass $\tau$ on ${\cali O}(M,\P)$,  such that 
$\lim_q \frac{1}{d^{kn_q}} \sum_{\gamma\in {\cali L}^r_{n_q}} \delta_\gamma =\cali M$ and $0<\tau \le 1$.
 
 We now fix $\lambda
 \in M$. By Briend-Duval theorem (\cite{BrDu}),   $\lim_n d^{-kn}\sum_{z\in {\cali R}_n(\lambda)} \delta_z =\mu_\lambda$ and, 
 setting $\sigma_\lambda:=e_{\lambda*}{\cali M}$, we have $\lim_q \frac{1}{d^{kn_q}} \sum_{z\in {\cali L}^r_{n_q}(\lambda)} \delta_z =\sigma_\lambda$. The positive measure $\sigma_\lambda$  is $f_\lambda$-invariant,  of mass $\tau$,  
and $\sigma_\lambda\le \mu_\lambda$. Writing $\mu_\lambda=\tau (\frac{\sigma_\lambda}{\tau }) + (1-\tau) \left(\frac{\mu_\lambda-\sigma_\lambda}{1-\tau }\right)$, we deduce from the ergodicity of $\mu_\lambda$ that $\sigma_\lambda =\tau \mu_\lambda$. Now, setting  ${\cali L}_n^r(\lambda):={\cali L}^r_n \cap (\{\lambda\}\times \P)$ and ${\cali R_n'(\lambda):=\cali R}_n(\lambda)\setminus 
 {\cali L}_n^r(\lambda)$, we get $\lim_q \frac{1}{d^{kn_q}} \sum_{z\in {\cali R}'_{n_q}(\lambda)} \delta_z =(1-\tau)\mu_\lambda$. 
 
Let $\ln_\epsilon$ be a family of smooth functions on $[0,+\infty[$ which converges pointwise to $\ln$ and such that $\ln_\epsilon \ge \ln$. Then
\begin{eqnarray*}
\tau \int_{\P} \ln_\epsilon \vert \Jac f_\lambda \vert\;\mu_\lambda =\int_{\P} \ln_\epsilon \vert \Jac f_\lambda'\vert\;\sigma_\lambda &=& \lim_q \frac{1}{d^{kn_q}} \sum_{z\in {\cali L}^r_{n_q}(\lambda)}  \ln_\epsilon \vert \Jac f_\lambda (z)\vert\\
&\ge&  \limsup_q \frac{1}{d^{kn_q}} \sum_{z\in {\cali L}^r_{n_q}(\lambda)}  \ln \vert \Jac f_\lambda (z)\vert
\end{eqnarray*}
and, making $\epsilon\to 0$,
\begin{eqnarray*}
\tau L(\lambda) \ge \limsup_q \frac{1}{d^{kn_q}} \sum_{z\in {\cali L}^r_{n_q}(\lambda)}  \ln \vert \Jac f_\lambda (z)\vert.
\end{eqnarray*}
Similarly, we have
\begin{eqnarray*}
(1-\tau) L(\lambda) \ge \limsup_q \frac{1}{d^{kn_q}} \sum_{z\in {\cali R}'_{n_q}(\lambda)}  \ln \vert \Jac f_\lambda (z)\vert.
\end{eqnarray*}
 On the other hand, $\lim_q  \frac{1}{d^{kn_q}} \sum_{z\in {\cali L}^r_{n_q}(\lambda)}  \ln \vert \Jac f_\lambda (z)\vert + \frac{1}{d^{kn_q}} \sum_{z\in {\cali R}'_{n_q}(\lambda)}  \ln \vert \Jac f_\lambda (z)\vert =L(\lambda)$  by the approximation formula. Thus,  $\tau L(\lambda) = \lim_q \frac{1}{d^{kn_q}} \sum_{z\in {\cali L}^r_{n_q}(\lambda)}  \ln \vert \Jac f_\lambda (z)\vert$,  which makes  $L$ appearing  as the pointwise limit of a locally uniformly bounded sequence of pluriharmonic functions. The function $L$ is therefore pluriharmonic on $M$.\qed\\

\end{document}